\documentstyle[11pt]{article}
\begin{document}
\newcommand{\p}{\parallel }
\makeatletter \makeatother
\newtheorem{th}{Theorem}[section]
\newtheorem{lem}{Lemma}[section]
\newtheorem{de}{Definition}[section]
\newtheorem{rem}{Remark}[section]
\newtheorem{cor}{Corollary}[section]
\renewcommand{\theequation}{\thesection.\arabic {equation}}

\title{{\bf Equivariant heat invariants of the Laplacian and nonmininmal operators on differential forms}}

\author{ Yong Wang\\
 }

\date{}
\maketitle

\begin{abstract}~~ In this paper, we compute the first two equivariant heat kernel coefficients of
the Bochner Laplacian on differential forms. The first two
equivariant heat kernel coefficients of the Bochner Laplacian with
torsion are also given. We also study the equivariant heat kernel
coefficients of nonmininmal operators on differential forms and get
the equivariant Gilkey-Branson-Fulling formula.\\
\\


\end{abstract}

\section{Introduction}
  \quad In [Do1], Donnelly computed heat kernel coefficients of the Bochner Laplacian with
torsion on differential forms. In [Do2], the equivariant heat kernel
asymptotics of the Laplacian on functions was established and the
first two equivariant heat kernel coefficients were evaluated.
Donnelly's results were used to study the asymptotic expansion of
the heat kernel for orbifolds in [DGGW]. On the other hand, in
[GBF], Gilkey, Branson and Fulling studied the heat equation
  asymptotics of nonmininmal operators on differential forms and
  computed the first two heat kernel coefficients. In [AV], the heat
  kernel expansion for a general nonmininmal operator on the spaces
  $C^{\infty}(\wedge^k)$ and $C^{\infty}(\wedge^{p,q})$ was studied. In [Gi],
  Gilkey studied the heat kernel asymptotics of nonmininmal operators
for manifolds with boundary. In [PC], the Gilkey-Branson-Fulling
formula was generalized to the case of the $h$-Laplacian. The
purpose of this paper is to compute the first two equivariant heat
kernel coefficients of the Bochner Laplacian and the Bochner
Laplacian with torsion on differential forms and give the
equivariant version of the Gilkey-Branson-Fulling formula. We
firstly compute the first two equivariant heat kernel coefficients
of the Bochner Laplacian on differential forms which are given in
Section 2. In section 3, we give the first two equivariant heat
kernel coefficients of the Bochner Laplacian with torsion on
differential forms. In Section 4.  we prove the equivariant
Gilkey-Branson-Fulling formula.\\

\section{The computation of the equivariant heat kernel
coefficients }
 \quad Let $M$ be a compact oriented Riemannian manifold of dimension $d$ and
 $T:M\rightarrow M$ be an isometry preserving orientation with fixed point set $\Omega$.
 $\Omega$ is the disjoint union of closed connected submanifolds $N$
 of dimension $n$. Let $TM$ and $T^*M$ be the tangent and cotangent
 bundles and let $\wedge^p$ be the bundle of exterior $p$-forms over
 $M$ and $\triangle$ be the Bochner Laplacion associated to the Levi-Civita connection on
 $\wedge^p$. $T$ induces a map $T^*$ on $\wedge^p$ which commutes
 with $\triangle$. If $\lambda$ is an eigenvalue of $\triangle$,
 then $T^*$ maps the $\lambda$ eigenspace into itself. We shall be
 interested in the sum $\sum_\lambda{\rm
 Tr}(T^*_\lambda)e^{-t\lambda}$, $t>0$ where Tr denotes the trace.
 Using the same discussions of Theorem 4.1 in [Do] (also see [Gi1]),
  we get\\

  \noindent {\bf Theorem 1} {\it  Let $T:M\rightarrow M$ be an isometry
  with fixed point set $\Omega$. Then there is an asymptotic
  expansion}
$$\sum_\lambda{\rm
 Tr}(T^*_\lambda)e^{-t\lambda}\sim \sum_{N\in\Omega}(4\pi
 t)^{-\frac{n}{2}}\sum_{k=0}^{\infty}t^k\int_Nb_k(T,a)d{\rm
 vol}_N(a),~~t\rightarrow 0, \eqno(2.1)$$
 {\it where the functions $b_k(T,a)$ depend only upon the germ of $T$
 and the Riemannian metric of $M$ near points $a\in N$.}\\

 Let $u_i(x,y)$ be the $i$-th term of the asymptotic
  expansion of the heat kernel of $e^{-t\bigtriangleup}$ and $A$
  denote the endomorphism induced by $T$ on the fibre of the normal
  bundle over $a$ and $B=(I-A)^{-1}$. Let $U_N$ be a sufficient
  small
  tubular neighborhood of $N$ and $\pi:U_N\rightarrow N$ be the
  projection. Let $x$ be coordinates for a normal coordinate chart
  on $\pi^{-1}(a)$ for $a\in N$ and $\overline{x}=x-T(x)$. Let $d{\rm vol}_M(x)=\psi(x)
  dx(\pi^*d{\rm vol}_N(a))$. By the Morse lemma, one can find a
  smooth coordinate change so that (see Appendix in [Do])
$$d^2(\overline{x}+T(x),T(x))=\sum_{i=1}^sy_i^2=|y|^2.$$
Let $|J(\overline{x},y)|$ denote the absolute value of the Jacobian
determinant of this change of variables. Let
$\Box_y=\sum_{i=1}^s\frac{\partial^2}{\partial y_i^2}$ and
$$h_i(x(y))=|{\rm det}B|{\rm
tr}[T^*u_i(T(x),x)]|J(\overline{x},y)|\psi(x),\eqno(2.2)$$ then by
the same discussions as in [Do], we have
$$b_k(T,a)=\sum_{j=0}^k\frac{1}{j!}\Box^j_y(h_{k-j})(0).\eqno(2.3)$$
Let $\tau(y,x):\wedge^p_x\rightarrow \wedge^p_y$ be the parallel
transport along the geodesic curve from $x$ to $y$. Fix $x$ and
suppose $y$ is in some normal coordinate neighborhood $w_j$ of $x$
and $g_{ij}=g(\partial/\partial w_i,\partial/\partial w_j)$. By
Chapter 2 in [BGV], we have $u_0(y,x)=({\rm det}
g_{ij})^{-\frac{1}{4}}\tau(y,x).$ By (2.3), $$b_0(T,a)=|{\rm
det}B|{\rm tr}[T^*u_0(T(a),a)]=|{\rm det}B|{\rm
tr}[T^*_a|_{\wedge^p}].$$ We write
$\widetilde{A}=\left(\begin{array}{cc}
\ I & 0 \\
 0 & A^t
\end{array}\right)$ and  define
$$\wedge^p\widetilde{A}(\theta_1\wedge\cdots\wedge\theta_p)=\widetilde{A}\theta_1\wedge\cdots\wedge
\widetilde{A}\theta_p,\eqno(2.4)$$ then
$${\rm tr}[\wedge^p\widetilde{A}]=\sum_{0\leq
p_1\leq p}\sum_{n+1\leq i_1<\cdots<i_{p_1}\leq
d}\sum_{l_1,\cdots,l_{p_1}}\varepsilon^{l_1,\cdots,l_{p_1}}_{i_1,\cdots,i_{p_1}}A_{i_1l_1}\cdots
A_{i_{p_1}l_{p_1}}\eqno(2.5), $$ where
$\varepsilon^{l_1,\cdots,l_{p_1}}_{i_1,\cdots,i_{p_1}}$ is the
generalized Kronecker symbol and $$b_0(T,a)=|{\rm det}B|{\rm
tr}[\wedge^p\widetilde{A}].\eqno(2.6)$$ By (2.3),
$$b_1(T,a)=|{\rm det}B|\left(h_1(0)+\Box_y({\rm
tr}[T^*u_0(T(x),x)]|J(\overline{x},y)|\psi(x))(0)\right).\eqno(2.7)$$
Denote $\tau_0=\sum_{a,b=1}^dR_{abab}$ and
$\rho_{ab}=\sum_{c=1}^dR_{acbc}$ the scalar curvature and Ricci
tensor of $M$. By Lemma 4.8.7 in [Gi1],
$$u_1(a,a)=\frac{\tau_0}{6},\eqno(2.8)$$
so
$$h_1(0)={\rm tr}[T^*_au_1(T(a),a)]=\frac{\tau_0}{6}{\rm
tr}[\wedge^p\widetilde{A}],\eqno(2.9)$$ By the Taylor expansions in
page 169 in [Do], we know that
$$\frac{\partial}{\partial y_i}|_{y=0}\left(({\rm det}
g_{ij})^{-\frac{1}{4}}|J(\overline{x},y)|\psi(x)\right)=0,\eqno(2.10)$$
so $$ \Box_y({\rm tr}[T^*u_0(T(x),x)]|J(\overline{x},y)|\psi(x))(0)
=\Box_y({\rm tr}[T^*\tau(T(x),x)])|_{y=0}$$ $$+{\rm
tr}[\wedge^p\widetilde{A}]\Box_y\left(({\rm det}
g_{ij})^{-\frac{1}{4}}|J(\overline{x},y)|\psi(x)\right)(0).\eqno(2.11)$$
In the following, we adopt the convention of summing Greek indices
$1\leq \alpha,\beta,\gamma\leq n$ from $1$ to $n$ and Latin indices
$n+1\leq i,j,k\leq d$ from $n+1$ to $d$. By the computations in
[Do,p.170], we know that
$$\Box_y\left(({\rm det}
g_{ij})^{-\frac{1}{4}}|J(\overline{x},y)|\psi(x)\right)(0)
=\frac{1}{6}\rho_{kk}+\frac{1}{3}R_{iksh}B_{ki}B_{hs}+\frac{1}{3}R_{ikth}B_{kt}B_{hi}-R_{k\alpha
h\alpha}B_{ks}B_{hs}.\eqno(2.12)$$ Now we compute ${\rm
tr}[T^*\tau(T(x),x)].$  Let $E=(E_1,\cdots,E_d)$ be an oriented
orthonormal frame field in a neighborhood of $a$ such that for $x\in
N$, $E_1(x),\cdots,E_n(x)$ are tangent to $N$ while the vector
fields $E_{n+1}(x),\cdots,E_d(x)$ are normal to $N$ and $E$ is
parallel along geodesics normal to $N$ and $dT$ is expressed as a
matrix-valued function $A(x)$ as $dTE(x)=E(Tx)A(x).$ Then
$A(a)=\left(\begin{array}{cc}
\ I & 0 \\
 0 & A
\end{array}\right)$. Let $x={\rm
exp}_a\left(\sum_{i=n+1}^dx_iE_i(a)\right).$ By Lemma 3.1 in [LYZ],
$A(x)=A(a)$. We consider $\wedge^p(T^*M)$ as the associated bundle
$SO(T^*M)\times_\mu \wedge^p({\bf R}^n)$ where the action $\mu$ is
defined by (2.4). Let $\sigma=(E_1^*,\cdots,E_d^*)$ be the local
section of $SO(T^*M)$ on $V$. Then a local section of
$\wedge^p(T^*M)$ on $V$ can be expressed as $[(\sigma,f)]$ where
$f:V\rightarrow \wedge^p({\bf R}^n)$ be a smooth function. Let
$$T^*[(\sigma(Tx),c)]=[(\sigma(x),\overline{T}^*(x)c)],\eqno(2.13)$$
where $\overline{T}^*(x):V\rightarrow {\rm End}({\bf R}^n)$. Let
$$\tau(Tx,x)[(\sigma(x),c)]=[(\sigma(Tx),\overline{\tau}^*(x)c)]. \eqno(2.14)$$
Then $\overline{T}^*(x)=\wedge^p\widetilde{A}$ is a constant
matrix-valued function. As in [LYZ, p.575], define the oriented
frame field $E^{Tx}$ over the patch $V$ by requiring that
$E^{Tx}(Tx)=E(Tx)$ and that $E^{Tx}$ be parallel along geodesics
through $Tx$ and a map $\Phi:V\rightarrow so(d)$ by
$E^{Tx}(x)=E(x)e^{\Phi(x)}.$ Then $E^{*,x}(Tx)=E^*(Tx)e^{-\Phi(x)}$,
that is
$$\tau(Tx,x)\sigma(x)=\sigma(Tx)e^{-\Phi(x)},\eqno(2.15)$$
so $\overline{\tau}^*(x)=\wedge^pe^{-\Phi(x)}.$ By Lemma 3.3 in
[LYZ], $\Phi=\left(\begin{array}{cc}
\ 0 & 0 \\
 0 & \Psi
\end{array}\right)$ and
$$\Psi_{ij}(x)=-\frac{1}{2}A_{rl}x_lx_sR_{rsij}(a)+O(|x|^3),$$
By $x_i=B_{ij}\overline{x_j}$ and $\overline{x_j}=y_j+O(|y|^3)$, so
$$\Psi_{ij}(x)=-\frac{1}{2}A_{rl}B_{lk}B_{sq}R_{rsij}(a)y_ky_q+O(|y|^3).$$
$$e^{-\Phi(x)}=1-\Phi(x)+O(|x|^3).$$
$${\rm
tr}[T^*\tau(T(x),x)]={\rm tr}[\wedge^p(W)]={\rm
tr}[\wedge^p(\widetilde{A}(1-\Phi(x)))]+O(|x|^3),\eqno(2.16)$$
$$C:=\Box_y({\rm tr}[\wedge^p(W)])|_{y=0}=\sum_{\delta=1}^{d-n}\sum_{1\leq
i_1<\cdots<i_p\leq
d}\sum_{l_1,\cdots,l_p}\varepsilon_{i_1,\cdots,i_p}^{l_1,\cdots,l_p}\sum_{s=1}^pW_{l_1i_1}\cdots
W_{l_{s-1}i_{s-1}}$$
$$\cdot(\frac{\partial^2}{\partial y^2_\delta}W_{l_{s}i_{s}})W_{l_{s+1}i_{s+1}}\cdots
W_{l_{p}i_{p}}|_{y=0}$$
$$=\sum_{\delta=1}^{d-n}\sum_{1\leq
i_1<\cdots<i_p\leq
d}\sum_{l_1,\cdots,l_p}\varepsilon_{i_1,\cdots,i_p}^{l_1,\cdots,l_p}\sum_{s=1}^p
A_{i_1l_1}\cdots A_{i_{s-1}l_{s-1}}A_{i_{s+1}l_{s+1}}\cdots
A_{i_{p}l_{p}}A_{rl}B_{l\delta}B_{v\delta}A_{ml_s}R_{rvmi_s}$$
$$=\sum_{\delta=1}^{d-n}\sum_{1\leq p_1\leq p}\sum_{n+1\leq
i_1<\cdots<i_{p_1}\leq
d}\sum_{l_1,\cdots,l_{p_1}}\varepsilon_{i_1,\cdots,i_{p_1}}^{l_1,\cdots,l_{p_1}}\sum_{s=1}^{p_1}
A_{i_1l_1}\cdots A_{i_{s-1}l_{s-1}}$$
$$\cdot A_{i_{s+1}l_{s+1}}\cdots
A_{i_{p_1}l_{p_1}}A_{rl}B_{l\delta}B_{v\delta}A_{ml_s}R_{rvmi_s}.\eqno(2.17)$$
So by (2.5),(2.7),(2.11),(2.12) and (2.17), we get
$$b_1(T,a)=|{\rm det }B|\left\{C+{\rm
tr}[\wedge^p\widetilde{A}](\frac{\tau_0}{6}+\frac{1}{6}\rho_{kk}\right.$$
$$\left.+\frac{1}{3}R_{iksh}B_{ki}B_{hs}+\frac{1}{3}R_{ikth}B_{kt}B_{hi}-R_{k\alpha
h\alpha}B_{ks}B_{hs})\right\}.\eqno(2.18)$$\\

\noindent{\bf Theorem 2} {\it The coefficient $b_k(T,a)$ of Theorem
1 is of the form $b_k(T,a)=|{\rm det }B|b'_k(T,a)$ where $b'_k(T,a)$
is an invariant polynomial in the components of $A$, $B$ and the
curvature tensor $R$ and its covariant derivative at $a$. In
particular, $b_0(T,a)$,~$b_1(T,a)$ are determined by (2.6) and
(2.18).}\\

\noindent {\bf Remark} Theorem 2 can be used to evaluate the heat
kernel coefficients of the Laplacian on forms on orbifolds as in
[DGGW].\\

\section{The computation of the equivariant heat kernel
coefficients of the Bochner Laplacian with torsion}
 \quad The Levi-Civita connection $\nabla$ is the unique torsion
 zero connection on $TM$ which preserves the metric. More generally,
 let $\overline{T}:TM\times TM\rightarrow TM$ be a skew-symmetric linear map,
 i.e. $\overline{T}(X,Y)=-\overline{T}(Y,X).$ Then there is exactly one metric-preserving
 connection $\overline{\nabla}=\nabla+Q$ on $TM$ with torsion tensor
 $\overline{T}$. Let
 $ Q(X,Y,Z)=g(Q(X,Y),Z);~$ $\overline{T}(X,Y,Z)=g(\overline{T}(X,Y),Z)$ and
 $Q_{ijk}=Q(E_i,E_j,E_k);~\overline{T}_{ijk}=\overline{T}(E_i,E_j,E_k),$ then
 $$Q_{ijk}=\frac{1}{2}(\overline{T}_{ijk}+\overline{T}_{kij}+\overline{T}_{kji}).\eqno(3.1)$$
Let $\overline{\nabla}^*={\nabla}^*+Q^*$ be the dual connection of
$\overline{\nabla}$ on $T^*M$ and
$\overline{Q}_{ijk}=g(Q^*(E_i)E^*_j,E^*_k)$, then
$\overline{Q}_{ikj}=-Q_{ijk}.$ $\overline{\nabla}^*$ induces a
connection on $\wedge^pT^*M$. We still denote it by
$\overline{\nabla}^*$. Let the Bochner Laplacian with torsion be
$$\overline{\triangle}=-\sum_{i=1}^d(\overline{\nabla}^*_{E_i}\overline{\nabla}^*_{E_i}
-\overline{\nabla}^*_{\nabla^L_{E_i}E_i}).$$ Let
$\overline{\tau}(y,x):\wedge^p_x\rightarrow \wedge^p_y$ be the
parallel transport along the geodesic curve from $x$ to $y$
associated to the connection $\overline{\nabla}^*$. Define the
oriented frame field $\overline{E}^{x,*}$ over the patch $V$ by
requiring that $\overline{E}^{x,*}(x)=E^*(x)$ and that
$\overline{E}^{x,*}$ be parallel along geodesics through $x$
associated to the connection $\overline{\nabla}^*$. Let $z$ be the
normal coordinates for the center $x$. Let
$$\overline{E}^{x,*}(z)=E^{x,*}(z)\overline{\Psi}(x,z);~\overline{\Psi}(x,Tx)=\overline{\Psi}(x).\eqno(3.2)$$
Nextly, we compute the Taylor expansion of $\overline{\Psi}(x).$ Let
$Q^*E^{x,*}(z)=E^{x,*}(z)A^x(z)$ and $L(x)=A^x(0)$, then
$$Q^*E^{*}(x)=E^{*}(x)L(x),~L(E_k)_{ij}=Q_{kji}.\eqno(3.3)$$
Let $\gamma$ be the geodesic curve from $x$ to $z$. By
$\nabla_{\dot{\gamma}}^{L,*}E^{x,*}=0,~\overline{\nabla}_{\dot{\gamma}}^{*}\overline{E}^{x,*}=0$
and (3.2), we have
$$\frac{d}{dt}\overline{\Psi}(x,{\gamma}(t))+A^x(\dot{\gamma}(t))\overline{\Psi}(x,{\gamma}(t))=0.\eqno(3.4)$$
By $\gamma(t)=tz;~\dot{\gamma}=\sum_{j=1}^d
z_j\frac{\partial}{\partial z_j},$
$$\sum_{j=1}^d
\frac{\partial\overline{\Psi}}{\partial
z_j}(x,tz)z_j+\sum_{j=1}^dz_jA^x(\frac{\partial}{\partial
z_j})(tz)\overline{\Psi}(x,tz)=0.\eqno(3.5)$$ By
$\overline{\Psi}(x,0)=id$ and setting $t=0$, we get
$$\sum_{j=1}^d
\frac{\partial\overline{\Psi}}{\partial
z_j}(x,0)z_j+\sum_{j=1}^dz_jA^x(\frac{\partial}{\partial
z_j})(0)=0.\eqno(3.6)$$ Taking the derivative about $t$ and setting
$t=0$, we get
$$\sum_{i,j=1}^d
\frac{\partial^2\overline{\Psi}}{\partial z_i\partial
z_j}(x,0)z_iz_j+\sum_{i,j=1}^dz_iz_jA^x(\frac{\partial}{\partial
z_j})(0)\frac{\partial\overline{\Psi}}{\partial
z_i}(x,0)+\sum_{i,j=1}^dz_iz_j\frac{\partial
A^x(\frac{\partial}{\partial z_j})}{\partial z_i}(0)=0.\eqno(3.7)$$
By (3.6) and (3.7), we have
$$\overline{\Psi}(x,z)=Id-\sum_{j=1}^dz_jA^x(\frac{\partial}{\partial
z_j})(0)+\frac{1}{2}\sum_{i,j=1}^dA^x(\frac{\partial}{\partial
z_j})(0)A^x(\frac{\partial}{\partial z_i})(0)z_iz_j$$
$$-\frac{1}{2}\sum_{i,j=1}^dz_iz_j\frac{\partial A^x(\frac{\partial}{\partial
z_j})}{\partial z_i}(0)+O(|z|^3).\eqno(3.8)$$ Let $Tx={\rm
exp}_x(\sum_{i=1}^du_iE_i(x))$ and $x=(a,c)$ be the orthogonal
coordinates in [LYZ,p.574], then by the proposition in [Yu,p.84], we
have
$$u_i=O(|c|^3),~1\leq i\leq n,~u_i=\overline{c_{i-n}}-c_{i-n}+O(|c|^3),~n+1\leq i\leq
d,\eqno(3.9)$$ where
$c=(\overline{c}_1,\cdots,\overline{c}_{d-n})A$. Then
$$\overline{\Psi}(x)=\overline{\Psi}(x,u)=
Id-\sum_{j=n+1}^du_jA^x(\frac{\partial}{\partial
z_j})(0)+\frac{1}{2}\sum_{i,j=n+1}^dA^x(\frac{\partial}{\partial
z_j})(0)A^x(\frac{\partial}{\partial z_i})(0)u_iu_j$$
$$-\frac{1}{2}\sum_{i,j=n+1}^du_iu_j\frac{\partial A^x(\frac{\partial}{\partial
z_j})}{\partial z_i}(0)+O(|u|^3).\eqno(3.10)$$ We know that
$$A^x(\frac{\partial}{\partial
z_j})(0)=A^x(0)(E_j(x))=L(x)(E_j(x)),\eqno(3.11)$$
$$\frac{\partial A^{(a,0)}(\frac{\partial}{\partial
z_j})}{\partial z_j}(0)=\frac{\partial
A^{(a,0)}(\frac{\partial}{\partial c_{j-n}})}{\partial
c_{i-n}}(0)=\frac{\partial L(E_j)}{c_{i-n}}(0),~n+1\leq i,j\leq d,
\eqno(3.12)$$
$$L(E_j)=L(E_j)(a,0)+\sum_{i=n+1}^d\frac{\partial
L(E_j)}{\partial c_{i-n}}(a,0)c_{i-n}+O(|c|^2),\eqno(3.13)$$ so
$$\overline{\Psi}(x)=Id-\sum_{j=n+1}^du_jL(E_j)(a,0)-\sum_{i,j=n+1}\frac{\partial
L(E_j)}{\partial c_{i-n}}(a,0)c_{i-n}u_j$$
$$+\frac{1}{2}\sum_{i,j=n+1}^dL(E_i)(a,0)L(E_j)(a,0)u_iu_j
-\frac{1}{2}\sum_{i,j=n+1}\frac{\partial L(E_j)}{\partial
c_{i-n}}(a,0)u_iu_j+O(|c|^3).\eqno(3.14)$$ By
$$u_i=-y_{i-n}+O(|y|^3),~c_{i-n}=B_{ij}y_{j-n}+O(|y|^3),~n+1\leq
i,j\leq d,\eqno(3.15)$$ we have
$$\overline{\Psi}(x)=Id+\sum_{j=n+1}^dy_{j-n}L(E_j)(a,0)+\sum_{i,j=n+1}\frac{\partial
L(E_j)}{\partial c_{i-n}}(a,0)y_{j-n}B_{ik}y_{k-n}$$
$$+\frac{1}{2}\sum_{i,j=n+1}^dL(E_i)(a,0)L(E_j)(a,0)y_{i-n}y_{j-n}
-\frac{1}{2}\sum_{i,j=n+1}\frac{\partial L(E_j)}{\partial
c_{i-n}}(a,0)y_{i-n}y_{j-n}+O(|y|^3).\eqno(3.16)$$ By (3.2),
$$\overline{E}^{x,*}(Tx)=E^{*}(Tx)e^{-\Phi(x)}\overline{\Psi}(x).\eqno(3.17)$$
Let
$$\overline{\tau}(Tx,x)[(\sigma(x),c)]=[(\sigma(Tx),\widetilde{\tau}^*(x)c)].
\eqno(3.18)$$ Let $\overline{b}_i(T,a)$ denote the equivariant heat
kernel coefficients of the Bochner Laplacian with torsion. Similar
to the discussions in Section 2, $\overline{b}_0(T,a)=b_0(T,a)$.
Similar to (2.18), we have
$$\overline{b}_1(T,a)=|{\rm det }B|\left\{\overline{C}+{\rm
tr}[\wedge^p\widetilde{A}](\frac{\tau_0}{6}+\frac{1}{6}\rho_{kk}\right.$$
$$\left.+\frac{1}{3}R_{iksh}B_{ki}B_{hs}+\frac{1}{3}R_{ikth}B_{kt}B_{hi}-R_{k\alpha
h\alpha}B_{ks}B_{hs})\right\}.\eqno(3.19)$$ where
$$\overline{C}=\Box_y({\rm tr}[T^*\widetilde{\tau}(T(x),x)])|_{y=0}
={\rm
tr}[\wedge^p(\widetilde{A}e^{-\Phi(x)}\overline{\Psi}(x))]:={\rm
tr}[\wedge^p\overline{W}].\eqno(3.20)$$ So
$$\overline{C}=\sum_{\delta=1}^{d-n}\sum_{1\leq
i_1<\cdots<i_p\leq
d}\sum_{l_1,\cdots,l_p}\varepsilon_{i_1,\cdots,i_p}^{l_1,\cdots,l_p}\frac{\partial^2}{\partial
y^2_\delta}[\overline{W}_{l_1i_1}\cdots
\overline{W}_{l_{p}i_{p}}]|_{y=0}$$
$$=\sum_{\delta=1}^{d-n}\sum_{1\leq i_1<\cdots<i_p\leq
d}\sum_{l_1,\cdots,l_p}\varepsilon_{i_1,\cdots,i_p}^{l_1,\cdots,l_p}
\left[2\sum_{1\leq m_1<m_2\leq p}\overline{W}_{l_1i_1}\cdots
\frac{\partial \overline{W}_{l_{m_1}i_{m_1}}}{\partial
y_\delta}\cdots \frac{\partial
\overline{W}_{l_{m_2}i_{m_2}}}{\partial
y_\delta}\cdots\overline{W}_{l_{p}i_{p}}\right.$$
$$\left.\sum_{1\leq m_3\leq p}\overline{W}_{l_1i_1}\cdots
\frac{\partial^2 \overline{W}_{l_{m_3}i_{m_3}}}{\partial
y^2_\delta}\cdots\overline{W}_{l_{p}i_{p}}\right]|_{y=0}$$ By (3.3)
and (3.16)
$$\frac{\partial\overline{\Psi_{\varepsilon
i_{m_1}}}}{\partial
y_\delta}|_{y=0}=-Q_{\delta+n~i_{m_1}\varepsilon};$$
$$\frac{\partial^2\overline{\Psi}_{a,b}}{\partial
y^2_\delta}|_{y=0}=2\sum_{j=n+1}^dQ_{\delta+n~ba,j}B_{j\delta+n}+\sum_{a_1=1}^dQ_{\delta+n~a_1a}Q_{\delta+n~ba_1}
-Q_{\delta+n~ba,\delta+n},\eqno(3.21)$$ then
$$\frac{\partial \overline{W}_{l_{m_1}i_{m_1}}}{\partial
y_\delta}|_{y=0}=-A_{\varepsilon
l_{m_1}}Q_{\delta+n~i_{m_1}\varepsilon};\eqno(3.22)$$
$$\frac{\partial^2 \overline{W}_{l_{m_3}i_{m_3}}}{\partial
y^2_\delta}|_{y=0}=A_{\varepsilon_1
l_{m_3}}[A_{rl}B_{l\delta}B_{v\delta}A_{m\varepsilon_1}R_{rvm
i_{m_3}}+2\sum_{j=n+1}^dQ_{\delta+n~i_{m_3}\varepsilon_1,j}B_{j\delta+n}$$
$$+\sum_{a_1=1}^d
Q_{\delta+n~a_1\varepsilon_1}Q_{\delta+n~i_{m_3}a_1}
-Q_{\delta+n~i_{m_3}\varepsilon_1,\delta+n}].\eqno(3.23)$$ So
$$\overline{C}=\sum_{\delta=1}^{d-n}\sum_{1\leq i_1<\cdots <i_p\leq
d}\sum_{l_1,\cdots,l_p}\varepsilon_{i_1,\cdots,i_p}^{l_1,\cdots,l_p}
\left[2\sum_{1\leq m_1<m_2\leq p}A_{l_1i_1}\cdots
\widehat{A}_{l_{m_1}i_{m_1}}\cdots\widehat{A}_{l_{m_2}i_{m_2}}
\cdots A_{l_{p}i_{p}}\right.$$
$$\cdot A_{\varepsilon_1
l_{m_1}}Q_{\delta+n~i_{m_1}\varepsilon_1}A_{\varepsilon_2
l_{m_2}}Q_{\delta+n~i_{m_2}\varepsilon_2}+\sum_{1\leq m_3\leq p}
A_{l_1i_1}\cdots \widehat{A}_{l_{m_3}i_{m_3}}\cdots A_{l_{p}i_{p}}$$
$$\cdot
A_{\varepsilon_1
l_{m_3}}[A_{rl}B_{l\delta}B_{v\delta}A_{m\varepsilon_1}R_{rvm
i_{m_3}}+2\sum_{j=n+1}^dQ_{\delta+n~i_{m_3}\varepsilon_1,j}B_{j\delta+n}$$
$$\left.+\sum_{a_1=1}^d
Q_{\delta+n~a_1\varepsilon_1}Q_{\delta+n~i_{m_3}a_1}
-Q_{\delta+n~i_{m_3}\varepsilon_1,\delta+n} \right]\eqno(3.24)$$

\noindent{\bf Theorem 3} {\it The coefficient $\overline{b}_k(T,a)$
is of the form $\overline{ b}_k(T,a)=|{\rm det
}B|\overline{b'}_k(T,a)$ where $\overline{b'}_k(T,a)$ is an
invariant polynomial in the components of $A$, $B$ and the curvature
tensor $R$ and the torsion tensor $\overline{T}$ and its covariant
derivative at $a$. In particular, $\overline{b}_1(T,a)$ are
determined by (3.24) and
(3.19).}\\

In the following, we define another Bochner's Laplacian with torsion
and compute its equivariant heat invariants. Let
$$\widehat{\triangle}=-\sum_{i=1}^d({\nabla}^*_{E_i}{\nabla}^*_{E_i}
-{\nabla}^*_{\overline{\nabla}_{E_i}E_i})=\triangle+Q_{iij}E_j+F.$$
Let $\widehat{\nabla}=\nabla-\frac{1}{2}(Q(E_i)E_i)^*$ be a
connection on $\wedge T^pM$ associated to a connection
$\widehat{\nabla}=\nabla-\frac{1}{2p}(Q(E_i)E_i)^*$ on $T^*M$. Then
$\widehat{\triangle}$ is a generalized Laplacian associated to the
connection $\widehat{\nabla}$. In this case,
$L(x)=-\frac{1}{2p}(Q(E_i)E_i)^*$ and $L(E_j)=-\frac{1}{2p}Q_{iij}$.
By Proposition 3.3, 3.5 and 4.3 [Do1], we have
$$\widehat{u}_1(x,x)=\frac{\tau_0}{6}-\frac{1}{2}T_{kjk,j}-\frac{1}{4}T_{kjk}T_{ljl},\eqno(3.25)$$
$$\widehat{b}_1(T,a)=|{\rm det }B|\left\{\widehat{C}+{\rm
tr}[\wedge^p\widetilde{A}](\frac{\tau_0}{6}-\frac{1}{2}T_{kjk,j}-\frac{1}{4}T_{kjk}T_{ljl}
+\frac{1}{6}\rho_{kk}\right.$$
$$\left.+\frac{1}{3}R_{iksh}B_{ki}B_{hs}+\frac{1}{3}R_{ikth}B_{kt}B_{hi}-R_{k\alpha
h\alpha}B_{ks}B_{hs})\right\}.\eqno(3.26)$$
$$\frac{\partial \overline{W}_{l_{m_1}i_{m_1}}}{\partial
y_\delta}|_{y=0}=\frac{1}{2p}A_{i_{m_1}l_{m_1}}Q_{ii\delta+n};\eqno(3.27)$$
$$\frac{\partial^2 \overline{W}_{l_{m_3}i_{m_3}}}{\partial
y^2_\delta}|_{y=0}=A_{\varepsilon_1
l_{m_3}}A_{rl}B_{l\delta}B_{v\delta}A_{m\varepsilon_1}R_{rvm
i_{m_3}}+A_{i_{m_3}
l_{m_3}}[-\frac{1}{p}\sum_{j=n+1}^dQ_{ii\delta+n,j}B_{j\delta+n}$$
$$+\frac{1}{4p^2}
Q_{ii\delta+n}Q_{kk\delta+n}
+\frac{1}{2p}Q_{ii\delta+n,\delta+n}].\eqno(3.28).$$ So
$$\widehat{C}=\sum_{\delta=1}^{d-n}\sum_{1\leq i_1<\cdots <i_p\leq
d}\sum_{l_1,\cdots,l_p}\varepsilon_{i_1,\cdots,i_p}^{l_1,\cdots,l_p}
\left[\frac{C_p^2}{2p^2}\sum_{1\leq m_1<m_2\leq p}A_{l_1i_1}\cdots
 A_{l_{p}i_{p}}Q_{ii\delta+n}Q_{kk\delta+n}
 \right.$$
$$+\sum_{1\leq m_3\leq p}
A_{l_1i_1}\cdots \widehat{A}_{l_{m_3}i_{m_3}}\cdots A_{l_{p}i_{p}}
\cdot A_{\varepsilon_1
l_{m_3}}[A_{rl}B_{l\delta}B_{v\delta}A_{m\varepsilon_1}R_{rvm
i_{m_3}}$$ $$+A_{l_1i_1}\cdots
 A_{l_{p}i_{p}}(-\sum_{j=n+1}^dQ_{ii\delta+n,j}B_{j\delta+n}
+\frac{1}{4p} Q_{ii\delta+n}Q_{kk\delta+n}
+\frac{1}{2}Q_{ii\delta+n,\delta+n})].\eqno(3.29).$$

\noindent{\bf Theorem 4} {\it The coefficient $\widehat{b}_k(T,a)$
is of the form $\widehat{ b}_k(T,a)=|{\rm det
}B|\widehat{b'}_k(T,a)$ where $\widehat{b'}_k(T,a)$ is an invariant
polynomial in the components of $A$, $B$ and the curvature tensor
$R$ and the torsion tensor $\overline{T}$ and its covariant
derivative at $a$. In particular, $\widehat{b}_1(T,a)$ are
determined by (3.26) and
(3.29).}\\

\section{The equivariant Gilkey-Branson-Fulling
formula}
 \quad Since $T$ is a preserving orientation isometry, then $T^*$
 commutes with $d$, $\delta$ and $\triangle$, so $T^*$ preserves the
 Hodge decomposition. Let
 $$\triangle^{(p)}_d=\triangle^{(p)}|_{{\rm Im}d};~\triangle^{(p)}_\delta=\triangle^{(p)}|_{{\rm
 Im}\delta}$$
and for each $t>0$
$$f_T(t,d^{(p)})={\rm Tr}(T^*e^{-t\triangle^{(p)}_d});~f_T(t,\delta^{(p)})={\rm
Tr}(T^*e^{-t\triangle^{(p)}_\delta}),$$
$$f_T(t,\triangle^{(p)})={\rm Tr}(T^*e^{-t\triangle^{(p)}});~
\beta_{T,p}={\rm Tr}(T^*|_{{\rm ker}(\triangle^{(p)})}).$$ Then we
have via the Hodge decomposition theorem:
$$f_T(t,\triangle^{(p)})=\beta_{T,p}+f_T(t,d^{(p)})+f_T(t,\delta^{(p)}).\eqno(4.1)$$
Similar to the nonequivariant case, we have
$$f_T(t,d^{(p)})=f_T(t,\delta^{(p-1)}).\eqno(4.2)$$
By (4.1) and (4.2), we have
$$f_T(t,\delta^{(p)})=\sum_{j\leq
p}(-1)^{p-j}[f_T(t,\triangle^{(j)})-\beta_{T,j}].$$ On the other
hand,
$$f_T(t,d^{(p)})=f_T(t,\delta^{(p-1)})=\sum_{j\leq
p-1}(-1)^{p-j}[f_T(t,\triangle^{(j)})-\beta_{T,j}].$$ Let
$D^{(p)}=a^2d\delta+b^2\delta d$ acting on $\wedge^p$ where $a\neq
0,~b\neq 0$. If we make the same computations for operator $D^{(p)}$
we obtain
$$f_T(t,D^{(p)})={\rm
Tr}(T^*e^{-tD^{(p)}})=\beta_{T,p}+f_T(a^2t,d^{(p)})+f_T(b^2t,\delta^{(p)})$$
$$=f_T(b^2t,\triangle^{(p)})+\sum_{j<p}(-1)^{p-j}[f_T(b^2t,\triangle^{(j)})
-f_T(a^2t,\triangle^{(j)})].\eqno(4.3)$$ By Lemma 1.8.2 in [Gi1], we
have
$$f_T(t,D^{(p)})=\sum_{N\in\Omega}(4\pi
 t)^{-\frac{n_N}{2}}\sum_{k=0}^{\infty}t^kb^T_{N,k}(D^{(p)})$$
$$f_T(b^2t,\triangle^{(j)})=\sum_{N\in\Omega}(4\pi
 t)^{-\frac{n_N}{2}}b^{-n_N}\sum_{k=0}^{\infty}t^kb^T_{N,k}(\triangle^{(j)})b^{2k}$$
$$f_T(a^2t,\triangle^{(j)})=\sum_{N\in\Omega}(4\pi
 t)^{-\frac{n_N}{2}}a^{-n_N}\sum_{k=0}^{\infty}t^kb^T_{N,k}(\triangle^{(j)})a^{2k}$$
By equating coefficients of $t^l$ in the asymptotic expansion in
(4.3), we
get\\

\noindent{\bf Theorem 5}
$$\sum_{N\in\Omega}(4\pi
 )^{-\frac{n_N}{2}}b^T_{N,l+\frac{n_N}{2}}(D^{(p)})=b^{2l}\sum_{N\in\Omega}(4\pi
 )^{-\frac{n_N}{2}}b^T_{N,l+\frac{n_N}{2}}(\triangle^{(p)})$$
 $$+\sum_{j<p}(-1)^{p-j}
(b^{2l}-a^{2l})\sum_{N\in\Omega}(4\pi
 )^{-\frac{n_N}{2}}b^T_{N,l+\frac{n_N}{2}}(\triangle^{(j)}).\eqno(4.4)$$

 \noindent {{\bf Remark 1.}} If $n_N={\rm constant}$ and taking
 $l=-\frac{n_N}{2}$, then we have
 $$\sum_{N\in\Omega}b^T_{N,0}(D^{(p)})=b^{-n}\sum_{N\in\Omega}b^T_{N,0}(\triangle^{(p)})
 +\sum_{j<p}(-1)^{p-j}(b^{-n}-a^{-n})\sum_{N\in\Omega}b^T_{N,0}(\triangle^{(j)}).\eqno4.5)$$
If $n_N={\rm constant}$ and taking
 $l=-\frac{n_N}{2}+1$, then we have
 $$\sum_{N\in\Omega}b^T_{N,1}(D^{(p)})=b^{-n+2}\sum_{N\in\Omega}b^T_{N,1}(\triangle^{(p)})
 +\sum_{j<p}(-1)^{p-j}(b^{-n+2}-a^{-n+2})\sum_{N\in\Omega}b^T_{N,1}(\triangle^{(j)}).\eqno(4.6)$$

 \noindent {{\bf Remark 2.}} Let $M$ be a K\"{a}hler manifold and $T$
 preserve the orientation and the canonical almost complex
 structure, then $T^*$ commutes with
 $\partial,\overline{\partial},\partial^*, \overline{\partial}^*$.
 Similar to Theorem 3, we can get the equivariant version of the
 expression of heat kernel coefficients of nonminimal operators on K\"{a}hler
 manifolds in [AV].\\

 \noindent{\bf Acknowledgement}~This work was supported by NSFC No.
10801027.\\

\noindent {\bf References}\\

\noindent[AV]S. Alexandrov, D. Vassilevich, {\it Heat kernel for
nonminimal operators on a Kahler manifold.} J. Math. Phys. 37
(1996), no. 11,
5715--5718.\\
\noindent[BGV] N. Berline; E. Getzler; M. Vergne, {\it Heat kernels
and Dirac operators.} Corrected reprint of the 1992 original.
Grundlehren Text Editions. Springer-Verlag, Berlin, 2004.\\
\noindent[Do1] H. Donnelly, {\it Heat equation asymptotics with
torsion,} Indiana Univ. Math. J. 34 (1985), no. 1, 105--113.\\
 \noindent[Do2] H. Donnelly, {\it Spectrum and the fixed point sets of isometries I.} Math. Ann. 224 (1976), no. 2,
 161--170.\\
\noindent[DGGW] E. Dryden; C. Gordon; S. Greenwald; D. Webb, {\it
Asymptotic expansion of the heat kernel for orbifolds.}
Michigan Math. J. 56 (2008), no. 1, 205--238.\\
 \noindent[LYZ] J. D. Lafferty, Y. L. Yu and W. P. Zhang, {\it A direct geometric proof of
Lefschetz fixed point formulas}, Trans.
AMS. 329(1992), 571-583.\\
 \noindent[G] P. Gilkey, {\it Heat content asymptotics of
nonminimal operators.} Topol. Methods Nonlinear Anal. 3 (1994), no. 1, 69--80.\\
\noindent[G1] P. Gilkey, {\it Invariance theory, the heat equation,
and the Atiyah-Singer index theorem.} Second edition. Studies in
Advanced
Mathematics. CRC Press, Boca Raton, FL, 1995.\\
 \noindent[GBF] P. Gilkey; T. Branson; S. Fulling, {\it Heat equation asymptotics of ``nonminimal'' operators
 on differential forms.} J. Math. Phys. 32 (1991), no. 8,
 2089--2091.\\
\noindent[PC] M. Puta; F. Cret, {\it A generalization of the
Gilkey-Branson-Fulling formula.} Proceedings of
  the Workshop on Global Analysis, Differential Geometry,
 Lie Algebras (Thessaloniki, 1997), 79--82, BSG Proc., 5, Geom. Balkan Press, Bucharest,
 2001.\\
\noindent[Yu] Y. Yu, {\it Trigonometry. II,} Acta Math. Sinica
(N.S.) 6 (1990), no. 1, 80--86.\\

 \indent School of Mathematics and
Statistics , Northeast Normal University, Changchun, Jilin 130024, China ;\\

 \indent E-mail: {\it wangy581@nenu.edu.cn}
\end {document}